\documentclass[12pt,a4paper]{amsart} 
\usepackage{epsf,amsfonts,amscd,latexsym}
\usepackage[T2A]{fontenc}
\usepackage[english]{babel} 
\def\hepsffile{\leavevmode\epsffile}

\title[Birational models of del Pezzo fibrations] 
{Birational models of del Pezzo fibrations}  
\thanks{This work was supported by the London Mathematical Society, 
grants RFBR no. 00--15--96085 and 02--01--00441, 
INTAS no. 00---0259 and 00--0269, and NWO-RFBR no. 047-008-005} 
\author{Mikhail Grinenko} 
\address{Steklov Mathematical Institute} 
\email{grin@mi.ras.ru} 
\date{}

\newtheorem{theorem}{\sc Theorem}[section]
\newtheorem{proposition}[theorem]{\sc Proposition}
\newtheorem{lemma}[theorem]{\sc Lemma}

\newtheorem{conjecture}[theorem]{\sc Conjecture}
\newtheorem{definition}[theorem]{\sc Definition}

\makeatletter

\@addtoreset{equation}{section} 
\newcommand{\l@abcd}[2]{\hbox to\textwidth{#1\dotfill #2}}

\makeatother

\newcommand*{\mybegintheorem}[1]{\begin{trivlist}\it%
  \item[\hspace{\labelsep}{\bf #1}]}
  \newcommand*{\myendtheorem}{\end{trivlist}}
  \newenvironment*{theorem*}{\mybegintheorem{Theorem.}}{\myendtheorem}
  \newenvironment*{proposition*}{\mybegintheorem{Proposition.}}{\myendtheorem}
  \newenvironment*{corollary*}{\mybegintheorem{Corollary.}}{\myendtheorem}
  \newenvironment*{definition*}{\mybegintheorem{Definition.}}{\myendtheorem}

\theoremstyle{remark} 
\newtheorem{remark}[theorem]{\sc Remark}
\newtheorem{example}[theorem]{\sc Example}

\renewcommand{\phi}{\varphi} 
\renewcommand{\epsilon}{\varepsilon}

\newcommand{\lra}{\longrightarrow}

\newcommand{\PN}{{{\mathbb P}^n}}
\newcommand{\PS}{{{\mathbb P}^6}}
\newcommand{\PF}{{{\mathbb P}^5}} 
\newcommand{\PQ}{{{\mathbb P}^4}}
\newcommand{\PT}{{{\mathbb P}^3}} 
\newcommand{\PTw}{{{\mathbb P}^2}}
\newcommand{\POn}{{{\mathbb P}^1}} 
\newcommand{\ZA}{{\mathbb Z}}
\newcommand{\QA}{{\mathbb Q}} 
\newcommand{\RA}{{\mathbb R}}
\newcommand{\CA}{{\mathbb C}} 
\newcommand{\PA}{{\mathbb P}}
\newcommand{\mc}{\mathcal}

\newcommand{\rk}{\mathop{\rm rk}\nolimits}
\newcommand{\Pic}{\mathop{\rm Pic}\nolimits}

\newcommand{\Bas}{\mathop{\rm Bas}\nolimits}

\newcommand{\cover}{\mathop{\lra}}

\newcommand{\Proj}{\mathop{{\bf Proj}\:}\nolimits}
\newcommand{\Spec}{\mathop{{\bf Spec}\:}\nolimits}

\newcommand{\eps}{\varepsilon}

\begin{document} 

\maketitle

\section{Preliminaries.}
\label{sec1}

Given an algebraic variety $X$, we can naturally attach some objects, 
e.g., its field of functions $k(X)$, the essential object in birational
geometry. So, assuming classification to be one of the most important 
problems in algebraic geometry, we may be asked to describe all 
algebraic varieties with the same field of functions, i.e., that are
birationally isomorphic to $X$. Of course, "all" is a too huge class, and 
usually we are restricted by projective and normal varieties (though 
non-proper or non-normal cases may naturally arise in some questions).
Typically there are two main tasks:

\begin{itemize}
\item[{\bf A.}] Given a variety $V$, one need to determine whether it 
is birational to another variety $W$.
\item[\bf B.] $V$ and $W$ are birational to each other, and one need
 to get a decomposition of a birational map between them into 
 "elementary links", i.e., birational maps that are simple enough.
\end{itemize}

The rationality problem is an essential example of task A. Recall that
a variety is said to be rational if it is birational to $\PN$ (or, 
which is the same, its field of functions is $k(x_1,\ldots,x_n)$).
As to task B, examples will be given below, let us only note that
varieties joined by "elementary links" should belong to a more or less
restrictive category (otherwise the task becomes meaningless), and
such a category have to meet the following requirement: any variety can 
be birationally transformed to one lying in it. Minimal models and Mori 
fibrations in dimension 3 give examples of such categories. Then, it
often happens that the indicated category, in its turn, is also too
large so as to be convenient, and we outline a subcategory of 
"good models" (e.g., relatively minimal surfaces). Here "good" means 
a class of varieties that are simple enough for describing, handling, 
classifying, and so on. Now let us look what we have in the first 
three dimensions.

\subsection{Curves.}
Normal algebraic curves are exactly smooth
ones. It is well known that a birational map between projective smooth
curves is an isomorphism, so the birational and biregular classifications
coincide in dimension 1. Projective spaces have no moduli, so $\POn$ is
a unique representative of rational curves.

\subsection{Surfaces.} 
It is very well known that any birational map between
smooth surfaces can be decomposed into a chain of blow-ups of points
and contractions of $(-1)$-curves without loss of the smoothness of
intermediate varieties. Now it is clear that we can successively contract 
(in any order) all $(-1)$-curves and get the so-called (relatively) 
minimal model (i.e., nothing to contract). So it is very convenient to 
use minimal models as the class of "good" models and the indicated blow-ups 
and contractions as elementary links. As to the rationality problem, 
the famous theorem of Castelnuovo says that a smooth surface $X$ is 
rational if and only if $H^1(X,{\mc O}_X)=H^0(X,2K_X)=0$. This is one of 
the most outstanding achievements of the classical algebraic geometry.

Summarizing results in dimensions 1 and 2, we can formulate the
rationality criterion as follows: $X$ is rational if and only if all
essential differential-geometric invariants, i.e., 
$H^0\left(X,(\Omega_X^p)^{\otimes m}\right)$, vanish.

\subsection{Threefolds.}
As soon as we get to dimension 3, the situation becomes much harder. Now
it isn't obvious at all what is "a good model". We may proceed with the
way as in dimension 2, i.e., contracting everything that can be contracted.
This is the viewpoint of Mori's theory. But then starting from a smooth 
variety, we may loose smoothness very quickly and even get a "very bad" 
variety with two Weil divisors intersecting by a point (this is the case 
of "a small contraction", i.e., birational morphism which is an isomorphism
in codimension 1). 

Nevertheless, Mori's theory states that there is the smallest category
of varieties which is stable under divisorial contractions and flips (the last
is exactly a tool which allows to "correct" small contractions). In what
follows we shall not need details of this theory, the reader can find them
in many monographs (e.g., \cite{KolMor}). We only point that $X$ belongs 
the Mori category if it is a projective normal variety with at most 
$\QA$-factorial terminal singularities. It means that every Weil divisor
is $\QA$-Cartier, i.e., it becomes Cartier as soon as we take it with
some multiplicity, and for every resolution of singularities
$\phi:Y\to X$ we have
$$
K_Y=\phi^*(K_X)+\sum a_iE_i,
$$
where all discrepancies $a_i$ are positive rational numbers, $E_i$
exceptional divisors, and "=" means "equal as $\QA$-divisors", i.e., 
multiplying by a suitable integer number, we get the linear equivalency. 
In particular, there exists the intersection theory on such varieties, 
which is very similar to the usual one, the difference is mostly that 
we must involve rational numbers as intersection indices.

The Mori category has some nice properties. In particular, the Kodaira
dimension is a birational invariant under maps in this category. We recall
that the Kodaira dimension $kod(X)$ is the largest dimension of images
under (rational) maps defined by linear systems $|mK_X|$. We are mostly
interested in studying varieties of negative Kodaira dimension, i.e.,
when $H^0(X,mK_X)=0$ for all $m>0$, because their birational geometry is
very non-trivial. 

From now on, we will only consider varieties of negative Kodaira dimension.
So, what are "minimal models" in the Mori theory for such varieties? These are
so-called Mori fiber spaces. By definition, a triple $\mu:X\to S$ is
a Mori fiber space if $X$ is projective, $\QA$-factorial and terminal,
$S$ is a projective normal variety with $\dim S<\dim X$, and $\mu$ is
an extremal contraction of fibering type, i.e., the relative Pickard
number $\rho(X/S)=\rk\Pic(X)-\rk\Pic(S)$ is equal to 1 and $(-K_X)$
is $\mu$-ample.

In dimension 3 (the highest dimension where the Mori theory has been proved)
we have three possible types of Mori fiber spaces (or, briefly, Mori 
fibrations) $\mu:X\to S$:
\begin{itemize}
\item[1)] $\QA$-Fano, if $\dim S=0$ (i.e., $S$ is a point);
\item[2)] del Pezzo fibration, if $\dim S=1$ (the fiber over the
generic point of $S$ is a del Pezzo surface of the corresponding degree);
\item[3)] conic bundle, if $\dim S=2$ (the fiber over the generic point
of $S$ is a plane conic).
\end{itemize}

In what follows, we often denote a Mori fibration $\rho:V\to S$ as
$V/S$ or even simply $V$, if it is clear which structure morphisms and
bases are meant.

\begin{figure}[htbp]
\begin{center}
\epsfxsize 10cm
\hepsffile{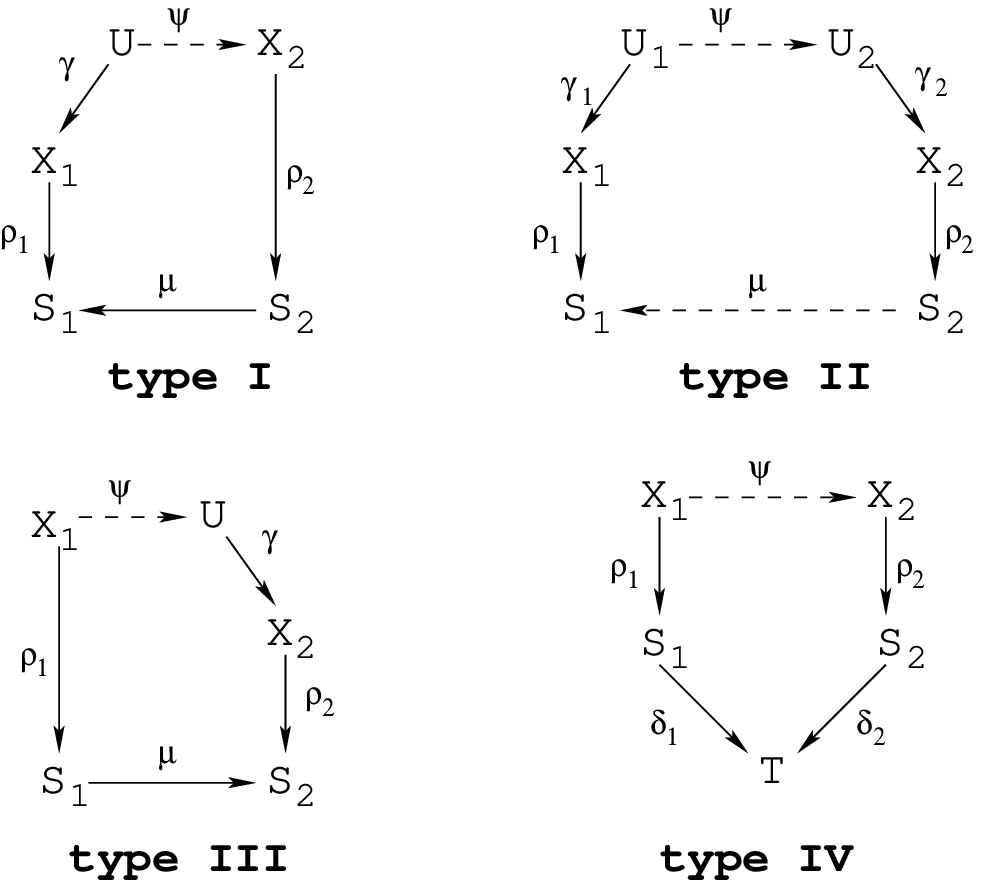}
\caption{}\label{london.fig}
\end{center}
\end{figure}

Factorization of birational maps between Mori fibrations is given by the 
Sarkisov program (it is proved in dimension 3, see \cite{Corti1}). The
essential assertion of this program is the following. Suppose we have
two Mori fibrations $V/S$ and $U/T$ and a birational map
$$
\begin{array}{ccc}
V & \stackrel{\chi}{\dasharrow} & U \\
\downarrow && \downarrow \\
S && T
\end{array}
$$ 
Then there exists a finite chain of birational maps
$$
\begin{array}{ccccccccccc}
X_0 & \stackrel{\chi_1}{\dasharrow} & X_1 &
\stackrel{\chi_2}{\dasharrow} & X_2 & \stackrel{\chi_3}{\dasharrow} &
\cdots & \stackrel{\chi_{N-1}}{\dasharrow} & X_{N-1} &
\stackrel{\chi_N}{\dasharrow} & X_N \\
\downarrow && \downarrow && \downarrow &&&& \downarrow &&\downarrow \\
S_0 && S_1 && S_2 &&&& S_{N-1} && S_N
\end{array}
$$
where $X_0/S_0,X_1/S_1,\ldots,X_N/S_N$ are Mori fibrations,
$X_0/S_0=V/S$, $X_N/S_N=U/T$, such that
$\chi=\chi_N\circ\chi_{N-1}\circ\ldots\circ\chi_2\circ\chi_1$ and any
of $\chi_i$ belongs to one of the four types of elementary links
listed below (see figure \ref{london.fig}). In all cases $X_1/S_1$ and
$X_2/S_2$ are Mori fibrations, $\psi$ is an isomorphism in codimension
1 (actually, a sequence of log-flips). Then, in the case of type I,
$\mu$ denotes a morphism with connected fibers, $\gamma$ is an
extremal divisorial contraction. Note that $\rho(S_1/S_0)=1$. In the
case of type II, $\gamma_1$ and $\gamma_2$ are extremal divisorial
contractions, $\mu$ is a birational map. Type III is inverse to type
I. Finally, $\delta_1$ and $\delta_2$ in type IV mean morphisms with
connected fibers, $T$ is a normal variety, and
$\rho(S_0/T)=\rho(S_1/T)=1$.

Let us note that in dimension 2 Mori fibrations are exactly relatively
minimal non-singular surfaces over points or curves, isomorphisms in
codimension 1 are actually biregular, and extremal divisorial contractions
are blow-downs of $(-1)$-curves. So all elementary links are very
simple. As to dimension 3, the links description just given is at
bottom all we know, so up to now the Sarkisov program mostly plays
a theoretical r\^ole.

The rationality problem for three-dimensional varieties also becomes 
enormously hard. During a long time, many mathematicians believed that 
it should be possible to find a rationality criterion more or less 
simple as in lower dimensions. But nearly simultaneously, in early 70th, 
three outstanding works of different authors which gave examples of 
unirational but non-rational varieties, appeared. These were Iskovskikh and 
Manin (\cite{IskMan}), Clemens and Griffiths (\cite{CG}), and Artin and 
Mumford (\cite{AM}). Recall that an algebraic variety is said to be 
unirational if there exists a rational map from projective space which 
is finite at the generic point. The matter is that all essential 
differential-geometric invariants vanish on unirational varieties as well,
so we can not even hope to find something like the rationality criterion
for curves and surfaces, combining these invariants. The reader may find
an excellent survey of the rationality problem in higher dimensions in
\cite{Isk1}.

Nevertheless, during the last 10 years, a considerable progress in the
birational classification problem in all dimension has been achieved, 
mostly due to the conception of birationally rigid varieties first 
formulated by A.V.Pukhlikov. The original version of the notion of 
birational rigidity is closely related to some technical features of the 
maximal singularities method (comparing of the canonical adjunction 
thresholds and all that), and still in active use. Another version
arised from the Mori fibrations theory. They are not identical but
coinside in many cases. In this paper we use the second one since it
is more understandable.

So, first we introduce the following useful notion: 
\begin{definition}
We say that Mori fibrations $V/S$ and $U/T$ ($\dim S=\dim T>0$) 
has the same Mori structure, if there exist birational maps 
$\chi:V\dasharrow U$ and $\psi:S\dasharrow T$ that make the following 
diagram to be commutative:
$$
\begin{array}{ccc}
V & \stackrel{\chi}{\dasharrow} & U \\
\downarrow && \downarrow \\
S & \stackrel{\psi}{\dasharrow} & T
\end{array}
$$
In $\QA$-Fano cases ($\dim S=\dim T=0$) we say that $V$ and $U$ has
the same Mori structure if they are biregular to each other; in other
words, for any birational map $\chi:V\dasharrow U$ there exists a
birational self-map $\mu\in Bir(V)$ such that the composition
$\chi\circ\mu$ is an isomorphism.
\end{definition}
In the sequel we will often use words "to be birational over the
base" instead of "to have the same Mori structure".

Now, since it is known that any threefold $X$ can be birationally
mapped onto a Mori fibration, one can formulate the classification
task as follows: to describe all Mori fibrations that are
birational to $X$ and not birational over the base to each other
(i.e., have different Mori structures). Clearly, varieties $X$ and $Y$
are birationally isomorphic, if and only if they have the same set of
Mori structures. By the way, note that only links of type II join Mori
fibrations with the same structure.

In the following cases the classification becomes especially easy:
\begin{definition}
A Mori fibration $V/S$ is said to be birationally rigid if it has a
unique Mori structure (i.e., any $U/T$ that is birational to $V$, is
birational over the base).
\end{definition}

Here are some examples of (birationally) rigid and non-rigid varieties.
But first let me note that any rigid variety is not rational. Indeed,
$\PT$ is birational to $\POn\times\PTw$, so we have at least three
different Mori structures: $\QA$-Fano ($\PT$ itself), del Pezzo
fibration $\POn\times\PTw\to\POn$, and conic bundle 
$\POn\times\PTw\to\PTw$. The simplest example of birationally rigid
$\QA$-Fano is a smooth quartic 3-fold in $\PQ$, this is nothing but
the result of Iskovskikh and Manin (\cite{IskMan}). Note that they
proved even more: any birational automorphism of a smooth quartic is 
actually biregular (and in general case, there are no non-trivial 
automorphisms at all). As to conic bundles, there exists the
following result of Sarkisov (\cite{Sark1}, \cite{Sark2}): any standard 
conic bundle $V\to S$ over a smooth rational surface $S$ with the 
discriminant curve $C$ is rigid if the linear system $|4K_S+C|$ is not empty
(i.e., $V/S$ has enough degenerations). Say, if 
$S$ is a plane, it is enough for $C$ to have degree 12 or higher.
Then, a smooth cubic 3-fold in $\PQ$ gives us an example of non-rigid
conic bundles. Indeed, the projection from a line lying on such a cubic,
realizes it as a conic bundle over a plane with a quintic as the
discriminant curve. 

Non-singular rigid and non-rigid del Pezzo fibrations will be
described in the main part of the paper. Almost all results were obtained
using the maximal singularities method (\cite{IskMan}, \cite{Pukh1}).
Here we shall not give any proves, the reader can find them in
\cite{Grin1}--\cite{Grin4}. 

Let me conclude this section with a little remark. Minimal models or
Mori fiber spaces are not convenient models all the way. In many cases,
there are more preferable classes of varieties. For examples, sometimes it
is useful to consider Gorenstein terminal varieties with numerically
effective anticanonical divisors (\cite{Al}, \cite{Corti2}).
In other words, we have to select each time.

\section{The rigidity problem for del Pezzo fibrations.}
\label{sec2}

Let $\rho:V\to S$ be a Mori fibration over, $S$ a smooth curve, $\eta$ the
generic point of $S$, and $V_{\eta}$ the fiber over the generic point.
So $V_{\eta}$ is a non-singular del Pezzo surface of degree
$d=(-K_{V_{\eta}})^2$ over the field of functions of $S$. Moreover,
this is a minimal surface since $\rho(V/S)=1$, so
$$
Pic(V_{\eta})\otimes\QA=\QA[-K_{V_{\eta}}].
$$
Suppose we have another Mori fibration $U/T$ and a birational map
$$
\begin{array}{ccc}
V & \stackrel{\chi}{\dasharrow} & U \\
\downarrow && \downarrow \\
S && T
\end{array}
$$
We would like to know as much as possible about this situation; in
particular, whether $\chi$ is birational over the base or not.

In first turn, we can assume $S$ to be rational. Indeed, suppose $S$
is a non-rational curve. Then $U$ can not be a $\QA$-Fano threefold.
The simplest reason is that the Pickard group of $V$ or any its 
resolution of singularities must contain a continuous part arising from
$S$, which is impossible for any birational model of $U$ (including $U$ 
itself). Or, we can use the fact that $U$ is rationally connected 
(see \cite{KolMM}), which is not true for $V$: there are no rational
curves lying across the fibers of $V/S$ (i.e., covering $S$).

Thus, $U/T$ is either a del Pezzo fibration or a conic bundle. Anyway,
the fibers of $U/T$ are covered by rational curves. Again, images of
these curves on $V$ can not cover $S$ and hence lie in the fibers of
$V/S$. So there exists a rational map from $T$ to $S$, which is
actually a morphism since the curve $S$ is smooth, and we have the
following commutative diagram:
$$
\begin{array}{ccc}
V & \stackrel{\chi}{\dasharrow} & U \\
\downarrow && \downarrow \\
S & \stackrel{\psi}{\longleftarrow} & T
\end{array}
$$
Denote $\zeta$ the generic point of $T$ and $U_{\zeta}$ the fiber over
the generic point. Assuming $U/T$ to be a del Pezzo fibration and
taking into account that $\chi$ is birational, we see easily that 
$\psi$ is an isomorphism, so it is possible to identify $\eta$ and
$\zeta$. Thus, $\chi$ induce the birational map $\chi_{\eta}$ between
two del Pezzo surfaces over $\eta$. Now suppose $U/T$ is a conic
bundle. Note that we may assume $\chi$ to be a link of type I.
The composition $U\to T\to S$ represents $U$ as a fibration over $S$,
so we can define the fiber $U_{\eta}$ over the generic point $\eta$.
It is clear that $U_{\eta}$ is a non-singular surface fibered on
conics with $Pic(U_{\eta})\simeq\ZA\oplus\ZA$ (see the
construction of the link). Moreover, we may assume $U_{\eta}$ to be
minimal. Indeed, $\rho(U_{\eta}/T_{\eta})=1$, so there are no
$(-1)$-curves in the fibers of $U_{\eta}\to T_{\eta}$. Hence any
$(-1)$-curve on $U_{\eta}$ has to be a section (which means that,
birationally, $U\to T$ is a $\POn$-bundle over $T$). Then simply
contract any of the $(-1)$-curves, and you obviously get a minimal 
del Pezzo surface, i.e., the previous case.
So, anyway, we have a birational map
$$
\chi_{\eta}: V_{\eta}\dasharrow U_{\eta}
$$
between two non-singular minimal surfaces defined over the field of
functions of $S$, and $U_{\eta}$ is either a del Pezzo surface or a
conic bundle. This situation was completely studied in V.A.Iskovskikh's paper
\cite{Isk2}, theorem 2.6. In particular, if $V_{\eta}$ is a del Pezzo
surface of small degree ($d\le 3$), then $U_{\eta}$ is always a del
Pezzo surface that is isomorphic to $V_{\eta}$, and $\chi_{\eta}$ is
either an isomorphism if $d=1$, or a composition of Bertini and/or 
Geiser involutions if $d=2,3$.

Thus, in what follows we assume $S=\POn$. Note that the field of
functions $\CA(S)$ of $S$ is the so-called $C_1$ field, and from
\cite{CT} it follows that $V_{\eta}$ always has a point over $\CA(S)$
(i.e., $V\to S$ has a section). By the way, $V_{\eta}$ has a Zariski
dense subset of such points (\cite{KolMM}). Yu.I.Manin proved 
(\cite{Man}) that del Pezzo surfaces of degree 5 or greater with a
point over a perfect field are always rational over the field itself.
This is exactly our case. So if $d\ge 5$, $V_{\eta}$ is rational over
$\CA(S)=\CA(\POn)$, hence $V$ is rational over $\CA$. Thus, $V$ is
birational to $\PT$; in particular, $V/S$ is always non-rigid.
Now suppose $d=4$. Then we can blow up  a section of $V\to S$ and,
after some log flips in the fibers, get a structure of conic bundle
(or, simply, blow up a point on $V_{\eta}$, as in theorem 2.6 of
\cite{Isk2}). Thus, the case $d=4$ is non-rigid too.

So we see that, from the viewpoint of the rigidity problem, the only
cases $S=\POn$, $d\le 3$, are interesting. The cases of degree 1 and 2
will be described in details in the rest of the paper. Here we only
outline some results in the case $d=3$, which is least studied.

\begin{theorem}[\cite{Pukh3}]
\label{dp3_th}
Let $V/\POn$ be a non-singular Mori fibration on cubic surfaces with
the simplest degenerations (the last means that all singular fibers
have at most one ordinary double point). If 1-cycles $N(-K_V)^2-f$ are 
non-effective for all $N>0$, where $f$ is the class of a line in a 
fiber, then $V/\POn$ is birationally rigid.
\end{theorem}
The indicated condition on 1-cycles (the so-called $K^2$-condition) is 
also sufficient for degree 1 and 2 (without assumptions on degenerations) 
but not necessary (see examples in \cite{Grin1}). The requirement on 
degenerations arised from the technical reasons and, most probably, can 
be omitted. The simplest examples of non-rigid fibrations on cubic 
surfaces are the following.

\begin{example}[\it Quartics with a plane] Let $X\in\PQ$ be a quartic
threefold containing a plane. It is easy to see that, in general case,
such a quartic contains exactly 9 ordinary double points lying on the
plane.  Note that $X$ is a Fano but not a $\QA$-Fano variety (i.e.,
not in the Mori category). Indeed, hyperplane sections through the
plane cut off residual cubic surfaces that intersect by these 9
singular points, which is impossible for $\QA$-factorial varieties. In
fact, the plane and any of these residual cubics generates the Weil
divisors group of $X$. Now blow up the plane as a subvariety in $\PQ$, 
and let $V$ be the strict transform of $X$. Then $V$ is nothing but a
Mori fibration on cubic surfaces over $\POn$, and the birational
morphism $V\to X$ is a small resolution with 9 lines lying over
the singular points of $X$. It is easy to check that we can produce
a flop simultaneously at these 9 lines, ant then contract the strict 
transform of the plane, which is $\PTw$ with the normal bundle 
isomorphic to ${\mc O}(-2)$. We get a $\QA$-Fano variety $U$, which
has the Fano index 1 and a singular point of index 2. This $U$ is the
complete intersection of two weighted cubic hypersurfaces in
$\PA(1,1,1,1,1,2)$. The described transformation $V/\POn\dasharrow U$
is nothing but a link of type 3 ($\psi$ is the flop, $\gamma$ is the
contraction of $\PTw$, see figure \ref{london.fig}). Conjecturally, 
$V/\POn$ and $U$ represent all Mori structures of $X$.
\end{example}

\begin{example}[\it Cubic threefolds]
Let $V\in\PQ$ be a smooth cubic hypersurface. It is well known that
$V$ is non-rational (\cite{CG}). Nevertheless, $V$ is very far from
rigid varieties. Indeed, the projection from any line on $V$ gives us 
a conic bundle over $\PTw$ with the discriminant curve of degree 5.
Then, $V$ is birational to a non-singular Fano variety of index 1 and
degree 14, which is a section of $G(1,5)\subset\PA^{14}$ by a linear
subspace of codimension 5. Finally, blow up any plane cubic
curve on $V$, and you get a structure of fibration on cubic surfaces.
So, $V$ has all possible types of Mori structures in dimension 3
(i.e., $\QA$-Fano, conic bundle, del Pezzo fibration).
\end{example}

\begin{example}[\it Hypersurfaces of bidegree $(m,3)$ in $\POn\times\PT$] 
Let $V=V_m\subset\POn\times\PT$ be a non-singular hypersurface of
bidegree $(m,3)$. From the Lefschetz hyperplane sections theorem it follows
that $\Pic(V)\simeq\ZA[-K_V]\oplus\ZA[F]$, where $F$ is the class of a
fiber, so $V/\POn$ is a Mori fibration on
cubic surfaces by means of the projection $\POn\times\PT\to\POn$.
Consider the second natural projection $\pi:\POn\times\PT\to\PT$. If
$m=1$ then, clearly, $\pi|_V$ gives us a birational map
$V\dasharrow\PT$, so $V$ is rational (and, thus, non-rigid) in this 
case. If $m\ge 3$, then for general $V$ the conditions of theorem
\ref{dp3_th} are satisfied, so $V/\POn$ is rigid. We consider the
most interesting case $m=2$. First, by the result of Bardelli
\cite{Bard}, general $V_m$'s are non-rational if $m\ge 2$. Now we show
that $V=V_2$ is non-rigid. Note that a general $V/\POn$ has exactly 27
sections that are fibers of the projection $\pi$. On the other hand,
let $t$ be any other fiber of $\pi$, then $t$ intersects $V$ at 2
points (not necessary different), and we can transpose these points. 
Thus, outside of the 27 sections, there exists an involution $\tau$,
so we have a birational self-map $\tau\in Bir(V)$. On the other hand,
we can simultaneously produce a flop $\chi$ centered at the indicated
sections and get another Mori fibration $U/\POn$, which is also a
hypersurface of bidegree $(2,3)$ in $\POn\times\PT$. This $\chi$ is a
link of type IV. Note that $\dim|3(-K_V)-F|=1$, this pencil has no
base components, and the strict transforms of its elements become the
fibers of $U/\POn$. The same is true for the corresponding pencil on
$U$. It only remains to notice that $\tau$ maps the pencil $|F|$ to
$|3(-K_V)-F|$, so we can consider $\chi$ as the Sarkisov resolution of
$\tau$. In \cite{Sob} it was shown that general $V_2$ have only these
two Mori structure.
\end{example}

We conclude this section with some remarks. First, in all
non-rigid cases, it is easy to check that the linear systems $|n(-K_V)-F|$
are non-empty and free from base components. In \cite{Grin1} the
following conjecture was formulated:
\begin{conjecture}
\label{main_conj}
Suppose $V/\POn$ is a non-singular Mori fibration on del Pezzo
surfaces of degree 1,2, or 3. Then $V/\POn$ is birationally rigid, if
and only if the linear systems $|n(-K_V)-F|$ are either empty or not free
from base components for all $n>0$.
\end{conjecture}

This conjecture was completely proved for degree 1 (\cite{Grin1}) and,
under some conditions of generality for several rigid cases, for
degree 2 (\cite{Grin1} and \cite{Grin2}). Moreover, for all cases that
are known to be rigid or non-rigid, the statement of the conjecture
holds. Finally, it does work in one direction even for a general 
Mori fibration on del Pezzo surfaces of degree 1:
\begin{theorem}[\cite{Grin3}, theorem 3.3]
If a Mori fibration $V/\POn$ on del Pezzo surfaces of degree 1 is
birationally rigid, then for all $n>0$ the linear systems
$|n(-K_V)-F|$ are either empty or not free from base components.
\end{theorem}

\section{Projective models of del Pezzo fibrations}
\label{sec3}

In this section we construct some simple projective models for
non-singular (actually, for Gorenstein) fibrations on del Pezzo
surfaces of degree 1 and 2. But first let me recall the simplest 
projective (and weighted projective) models of del Pezzo surfaces of 
degree 1 and 2. 

\subsection{Models of del Pezzo surfaces of degree 1 and 2}
Consider first the degree 1. So, let $X$ be a non-singular del Pezzo 
surface of degree $d=K_X^2=1$. It is easy to see that a general member 
of the linear system $|-K_X|$ is non-singular, all its elements are
irreducible and reduced curves of genus 1, and $|-K_X|$ has a unique
base point $P$. Choose a non-singular curve $C\in|-K_X|$. Using the
exact sequence
$$
0\lra{\mc O}_X\left(-(i-1)K_X\right)\lra{\mc O}_X\left(-iK_X\right)
\lra{\mc O}_C(-iK_X)\lra 0
$$ 
and the Kodaira vanishing theorem, we see that 
$$
H^0\left(X,-iK_X\right)\lra H^0\left(C,(-iK_X)|_C\right)\lra 0
$$
are exact for all $i\ge 0$. So we have a surjective map of the
graded algebras
$$
{\mc A}_X=\bigoplus_{i\ge 0}H^0\left(X,-iK_X\right)\lra
{\mc A}_C=\bigoplus_{i\ge 0}H^0\left(C,(-iK_X)|_C\right)
$$
which preserves the graduation. Clearly, ${\mc A}_C$ is generated by
elements of degree $r\le 3$ (note that $\deg (-3K_X)|_C=3$, i.e.,
$(-3K_X)|_C$ is very ample), so ${\mc A}_X$ is generated by elements
of degree not higher than 3, too. Taking into account that
$$
\begin{array}{ccc}
h^0(X,-K_X)=2, & h^0(X, -2K_X)=4, & h^0(X,-3K_X)=7,
\end{array}
$$
we can write down the generators as follows:
$$
\begin{array}{c}
H^0(X,-K_X)=\CA<x,y>, \\
H^0(X,-2K_X)=\CA<x^2,xy,y^2,z>, \\
H^0(X,-3K_X)=\CA<x^3,x^2y,xy^2,y^3,zx,zy,w>,
\end{array}
$$
where we agree for $x$ and $y$ to be elements of weight 1, $z$ of
weight 2, and $w$ of weight 3. Thus the homogeneous components of the
graded algebra ${\mc A}_X$ are generated by monomials $x^iy^jz^kw^l$
of the corresponding degree. Now it is clear that $X$ can be embedded
as a surface in the weighted projective space $\PA(1,1,2,3)$:
$$
X=\Proj{\mc A}_X=\Proj\CA[x,y,z,w]/I_X\subset=\PA(1,1,2,3),
$$
where $I_X$ is a principal ideal generated by a homogeneous element of
degree 6. Indeed, $h^0(X,-6K_X)=22$, but the dimension of the
homogeneous component of $\CA[x,y,z,w]$ of degree 6 is equal to 23, so
there exists exactly one linear relation in it. This relation is
nothing but an equation of $X$. It only remains to
notice that $X$ necessarily avoids singular points of $\PA(1,1,2,3)$,
whence we can put down the equation of $X$ as follows (stretching the
coordinates if needed):
$$
X=\{w^2+z^3+zf_4(x,y)+f_6(x,y)=0\}\subset\PA(1,1,2,3),
$$
where $f_i$ are homogeneous polynomials of the corresponding degree.

It was the way to embed $X$ into a very well-known threefold as a
surface. However, there is another useful model of $X$. Let
$$
\PA(1,1,2,3)\dasharrow\PA(1,1,2)
$$
be a (weighted) projection from the point $(0:0:0:1)$. This point does
not lie on $X$, and as it follows from the equation of $X$, 
the restriction of the projection on $X$ gives us a morphism of degree 2
$$
\phi:X\to\PA(1,1,2).
$$
We can naturally identify $\PA(1,1,2)$ with a non-degenerated quadratic
cone $Q\subset\PT$. The base point of $|-K_X|$ lies exactly over the
vertex of the cone. The ramification divisor of $\phi$ is defined by
an equation of degree 3 in $Q=\PA(1,1,2)$. We can also describe this
picture as follows: $|-2K_X|$ has no base points and defines a morphism, 
which is exactly our $\phi$. Thus, we can obtain $X$ as a double cover of a
quadratic cone in $\PT$ branched along a non-singular cubic section that 
does not contain the vertex of the cone.

Now let us turn back to a non-singular Mori fibration $V/\POn$ on del
Pezzo surfaces of degree 1. Clearly, a general fiber of $V$ is smooth.
But what if $V/\POn$ has degenerations? Moreover, as we will see
later, it must have degenerations by the relative Pickard number
arguments. So suppose $X$ be a singular fiber of $V$. Then we can see
that $X$ is a projective Gorenstein (i.e., the dualizing sheaf is
invertible, as it follows from the adjunction formula) Cohen-Macaulay
irreducible reduced (since $V$ is non-singular) 2-dimensional scheme,
$K_X$ is anti-ample and $K_X^2=1$.

First suppose $X$ is normal. The situation was completely studied in 
\cite{HiWa}, and I can simply refer the reader to it. But in order to
deal with projective models of $X$, we only need the fact that $|-K_X|$ 
contains a non-singular elliptic curve (\cite{HiWa}, proposition 4.2), 
and then one can repeat arguments of the non-singular case word by
word. So $X$ can be defined by a homogeneous equation of degree 6 in
$\PA(1,1,2,3)$ (and as before, $X$ avoids singular points of the
weighted projective space since they are non-Gorenstein), or as a
double cover of a quadratic cone in $\PT$, branched along a (now
singular) cubic section not passing through the vertex of the cone. 
As to this cubic section, it may be reducible, but all its components 
must be reduced. It only remains to add that the base point of
$|-K_X|$ is always non-singular, and $X$ has either only Du
Val singularities, or a minimally elliptic singular point (in local
coordinates defined by $r^2+p^3+q^6=0$).

The non-normal case is much harder. Nevertheless, in \cite{Reid} it
was shown that the projective properties of the anticanonical linear
systems on $X$ are the same as in the non-singular case. In
particular, the base point of $|-K_X|$ is non-singular, $X$ can be
embedded into $\PA(1,1,2,3)$ exactly in the same way, or twice cover a
quadratic cone. The unique difference is that the cubic section has to
contain a non-reduced component. 

We will summarize all needed results concerning to the case $d=1$ 
a bit later, now let us consider del Pezzo surfaces of degree 2. This
case is fairly similar to the previous one, so we only outline the key
points. Consider a non-singular del Pezzo surface $X$, $K_X^2=2$. Then
$|-K_X|$ is base points free and contains a non-singular element $C$. As
before, we have a surjective map ${\mc A}_X\to{\mc A}_C$, but now
these algebras are generated by elements of degree not higher than 2.
We see that $h^0(X,-K_X)=3$ and $h^0(X,-2K_X)=7$, so suppose
$$
\begin{array}{c}
H^0\left(X,-K_X\right)=\CA<x,y,z>, \\
H^0\left(X,-2K_X\right)=\CA<x^2,y^2,z^2,xy,xz,yz,w>,
\end{array}
$$
where $x$, $y$, and $z$ are of weight 1, $w$ is of weight 2. Then,
$$
X=\Proj{\mc A}_X=\Proj\CA[x,y,z,w]/I_X\subset\PA(1,1,1,2),
$$
the principal ideal $I_X$ is generated by a homogeneous element of
degree 4: indeed, $h^0(X,-4K_X)=21$, but the dimension of the homogeneous
component of $\CA[x,y,z,w]$ of degree 4 is equal to 22. Since $X$
avoids the singular point of $\PA(1,1,1,2)$, we get the following
equation of $X$:
$$
X=\{w^2+f_4(x,y,z)=0\}\subset\PA(1,1,1,2),
$$
where $f_4$ is a homogeneous polynomial of degree 4. The projection
$\PA(1,1,1,2)\dasharrow\PA(1,1,1)=\PTw$ from the point $(0:0:0:1)$,
being restricted to $X$, gives a double cover of $\PTw$ branched over
a non-singular quartic curve. This morphism can be also defined by the 
linear system $|-K_X|$.

Projective models of singular (normal and non-normal) Gorenstein del
Pezzo surfaces of degree 2 have the same construction as in the
non-singular case: they are defined by a quartic equation in
$\PA(1,1,1,2)$ and do not contain the point $(0:0:0:1)$. Or, we can
construct them as double coverings of $\PTw$ branched over quartic
plane curves, but now these quartics have singular points and/or
non-reduced components (\cite{HiWa}, \cite{Reid}). Briefly speaking, 
the situation is the same as in the case $d=1$. In particular, normal
surfaces of this type may have either only Du Val singularities, or a
unique minimally elliptic singularity locally defined by an equation
$r^2+p^4+q^4=0$.

We summarize the results about projective models of del Pezzo surfaces
of degree 1 and 2 as follows:

\begin{proposition}
\label{dp12surf_prop}
Let $X$ be a projective Gorenstein irreducible reduced del Pezzo
surface of degree 1 or 2. 

Then, in the case of degree 1, $|-K_X|$ has a unique base point, which
is non-singular on $X$, $|-2K_X|$ is base points free, 
and $|-3K_X|$ is very ample and embeds $X$ into $\PA^6$ as a surface
of degree 9. The suitable choice of the coordinates $[x,y,z,w]$ of 
weights $(1,1,2,3)$ in $\PA(1,1,2,3)$ allows to define $X$ by an
equation of degree 6:
$$
w^2+z^3+zf_4(x,y)+f_6(x,y)=0,
$$
where $f_i$ are homogeneous polynomials of degree $i$. $|-2K_X|$
defines a finite morphism $\phi$ of degree 2:
$$
\phi: X \cover_{R_Q}^{2:1} Q\subset\PT,
$$
where $Q$ is a non-degenerated quadratic cone, $R_Q=R|_Q$ is the
ramification divisor, $R\subset\PT$ is  a cubic that does not pass
through the vertex of $Q$.

In the case of degree 2, $|-K_X|$ is base points free, $|-2K_X|$ is
very ample and embeds $X$ into $\PA^6$ as a surface of degree 8. $X$
can be defined as a surface by an equation of degree 4:
$$
w^2+f_4(x,y,z)=0,
$$
where $f_4$ is homogeneous of degree 4, $[x,y,z,w]$ are the
coordinates of weights $(1,1,1,2)$ in $\PA(1,1,1,2)$. Finally, $|-K_X|$
defines a finite morphism of degree 2
$$
\phi: X \cover_{R}^{2:1} \PTw,
$$
where $R$ is the ramification curve, $\deg_{\PTw}R=4$. 
\end{proposition}

Now it is clear that models of fibrations on del Pezzo surfaces over
$\POn$ can be obtained as the relative version of the constructions
just introduced.

\subsection{Models of fibrations on del Pezzo surfaces of degree 1}
\label{subsec32}
Let $\rho:V\to\POn$ be a non-singular fibration on del Pezzo surfaces
of degree 1. Consider a sheaf of graded ${\mc O}_{\POn}$-algebras
$$
{\mc A}_V=\bigoplus_{i\ge 0}\rho_*{\mc O}_V(-iK_V).
$$
Since $(-K_V)$ is $\rho$-ample (i.e., on each fiber its restriction
gives anticanonical divisor, which is ample by the assumption), then
clearly
$$
   V=\Proj_{\POn}{\mc A}_V.
$$
Comparing the situation with proposition \ref{dp12surf_prop}, we see that 
there exist an algebra ${\mc O}_{\POn}[x,y,z,w]$ of polynomials over 
${\mc O}_{\POn}$ graded according to the weights $(1,1,2,3)$, and a
sheaf of principal ideals ${\mc I}_V$ in it such that
$$
V=\Proj_{\POn}{\mc A}_V=\Proj_{\POn}{\mc O}_{\POn}[x,y,z,w]/{\mc I}_V
\subset\PA_{\POn}(1,1,2,3),
$$
where $\PA_{\POn}(1,1,2,3)=\Proj_{\POn}{\mc O}_{\POn}[x,y,z,w]$. In
other words, $V/\POn$ is defined by a (weighted) homogeneous
polynomial (with coefficients in ${\mc O}_{\POn}$) in the
corresponding weighted projective space over $\POn$. Note that
$\PA_{\POn}(1,1,2,3)$ has two distinguished sections along which it is
singular, so $V$ does not intersect these sections.

Now let us construct the second model of $V/\POn$. All details and
proofs can be found in \cite{Grin1}, section 2.

In what follows, we assume that
$$
Pic(V)=\ZA[-K_V]\oplus\ZA[F],
$$
where $F$ denotes the class of a fiber. Notice that $V/\POn$ has a
distinguished section $s_b$ that intersects each fiber at the base
point of the anticanonical linear system. We can also define it as
$$
s_b=\Bas |-K_V+lF|
$$
for all $l\gg 0$. Then, $\rho_*(-2K_V+mF)$ is a vector bundle of rank
4 over $\POn$. We can choose m such that
$$
{\mc E}=\rho_*(-2K_V+mF)\simeq 
{\mc O}\oplus{\mc O}(n_1)\oplus{\mc O}(n_2)\oplus{\mc O}(n_3)
$$
for some $0\le n_1\le n_2\le n_3$. Denote 
$$
b=n_1+n_2+n_3.
$$
Let $X=\Proj {\mc E}$ be the corresponding $\PT$-fibration over
$\POn$, $\pi:X\to\POn$ the natural projection. Denote $M$ the class of
the tautological bundle on $X$ (i.e., $\pi_*{\mc O}(M)={\mc E}$), $L$
the class of a fiber of $\pi$, $t_0$ the class of a section that
corresponds to the surjection ${\mc E}\to{\mc O}\to 0$, and $l$ the
class of a line in a fiber of $\pi$. Note that $t_0$ is nothing but a 
minimally twisted over the base effective irreducible section of
$\pi$. Also $t_0$ can be obtained from the conditions $M\circ t_0=0$,
$L\circ t_0=1$. So, we have described all generators of the groups of
1- and 3-dimensional cycles on $X$.

Again, looking at proposition \ref{dp12surf_prop}, we see that there
exist a threefold $Q\subset X$ fibered into quadratic cones (without
degenerations) with a section $t_b$ as the line of the cone vertices, 
a divisor $R$ fibered into cubic surfaces and such that
its restriction $R_Q=R|_Q$ does not intersect $t_b$, and a finite 
morphism $\phi:V\to Q$ of degree 2 branched over $R_Q$, such that the
following diagram becomes commutative:
$$
\begin{CD}
V @>{\phi}>{R_Q,2:1}> Q \subset X\\
@V{\rho}VV @VV{\pi=\pi|_Q}V \\
\POn @= \POn
\end{CD}
$$
In order to make precise the construction, suppose $t_b\sim t_0+\eps
l$. Obviously, $\eps$ is some non-negative integer number.
\begin{lemma}[\cite{Grin1}, lemma 2.2]
The only following cases are possible:
\begin{itemize}
\item[1)] $\eps=0$, i.e., $t_b=t_0$, and then $2n_2=n_1+n_3$, 
 $n_1$ and $n_3$ are even, $Q\sim 2M-2n_2L$, and $R\sim 3M$;
\item[2)] $\eps=n_1>0$, and then $n_3=2n_2$, $n_1$ is even,
       $n_2\ge 3n_1$, $Q\sim 2M-2n_2L$, and $R\sim 3M-3n_1L$.
\end{itemize}
\end{lemma}
This lemma shows that the numbers $n_1$, $n_2$, and $n_3$ completely 
define the construction of $V/\POn$ (but up to moduli, of course).
Note that these numbers are not free from relations.

By the way, we can always assume $b>0$, i.e., $n_3>0$. Indeed,
otherwise $Q\sim 2M$ and $\eps=0$, so $V$ is nothing but the direct
product of $\POn$ and a del Pezzo surface, hence $V$ is rational and,
moreover, not a Mori fibration, because $\rho(V/\POn)$ coincide with
the Pickard number of the del Pezzo surface in this case and thus
greater than 1.

It only remains to consider some formulae and relations that allow to
identify $V/\POn$ more or less easy. First, let us note that a surface
$G\sim (M-n_2L)\circ (M-n_3L)$ must lie in $Q$. This is a minimally
twisted over the the base ruled surface on $Q$. Denote $G_V=\phi^*(G)$
its pre-image on $V$. Geometrically, $G_V$ is a minimally twisted
fibration on curves of genus 1 on $V$. Then, $t_0$ lies always on $Q$,
the fibers of $Q$ contain lines of the class $l$, so the classes
$s_0=\frac12\phi^*(t_0)$ and $f=\frac12\phi^*(l)$ are well defined.
Note that $f$ is the class of the anticanonical curves in the fibers
of $V$, and $s_0$ always has an effective representative in 1-cycles
on $V$. By the way, it is easy to see that the Mori cone is generated
by $s_0$ and $f$:
$$
\overline{\bf NE}(V)=\RA_+[s_0]\oplus\RA_+[f].
$$

It is not very difficult to compute the normal bundle of $s_b$:
\begin{itemize}
\item[] ${\mc N}_{s_b|V}\simeq{\mc O}(-\frac12n_1)\oplus{\mc O}
(-\frac12n_3)$, if $\eps=0$,
\item[]$ {\mc N}_{s_b|V}\simeq{\mc O}(n_1-\frac12n_3)\oplus{\mc O}(n_1)$,
  if $\eps=n_1>0$.
  \end{itemize}
So, in fact, we see that $V/\POn$ is completely defined (again, up to 
moduli) by ${\mc N}_{s_b|V}$.

The following table contains the essential information about divisors
and intersection indices on $V$:

\medskip
{
\renewcommand{\arraystretch}{1.5}
\begin{tabular}{|c|c|}
\hline
$\eps=0$ & $\eps=n_1>0$ \\
\hline
$n_1+n_3=2n_2$ & $n_3=2n_2$, $n_2\ge 3n_1$ \\
\hline
$s_b\sim s_0$ & $s_b\sim s_0+\frac12f$ \\
\hline
${\mc N}_{s_b|V}\simeq{\mc O}(-\frac12n_1)\oplus{\mc O}(-\frac12n_3)$
 & ${\mc N}_{s_b|V}\simeq{\mc O}(n_1-\frac12n_3)\oplus{\mc O}(n_1)$ \\
\hline
$K_V=-G_V+(\frac12n_1-2)F$ & $K_V=-G_V-(\frac12n_1+2)F$ \\
\hline
$K_V^2=s_0+(4-n_2)f$ & $K_V^2=s_0+(4+\frac32n_1-n_2)f$ \\
\hline
$s_0\circ G_V=-\frac12n_3$ & $s_0\circ G_V=-\frac12n_3$ \\
\hline
$(-K_V)^3=6-2n_2$ & $(-K_V)^3=6+2n_1-2n_2$ \\
\hline
\end{tabular}
}

\begin{remark}
\label{Gorenst_rem1}
The arguments still work if we assume $V$ to be only Gorenstein, not
necessary non-singular. The only point is that we must require all
fibers to be reduced. Indeed, all we need is the statement of
proposition \ref{dp12surf_prop}. It easy to see that in the Gorenstein
case all fibers of $V/\POn$ are irreducible (and reduced by the
requirement), so we are still under the conditions even if $V$ is 
singular. In particular, we have the same projective models as above, 
the same formulae, and so on. 
\end{remark}

\subsection{Models of fibrations on del Pezzo surfaces of degree 2}
\label{subsec33}
Let $\rho:V\to\POn$ be a non-singular fibration on del Pezzo surfaces
of degree 2. As in the case $d=1$, consider the direct image of the
anticanonical algebra
$$
{\mc A}_V=\bigoplus_{i\ge 0}\rho_*{\mc O}_V(-iK_V).
$$
This is a sheaf of graded algebras over ${\mc O}_{\POn}$, and
$V=\Proj_{\POn}{\mc A}_V$. As before, we see that there exist
a polynomial algebra ${\mc O}_{\POn}[x,y,z,w]$ graded with the weights
$(1,1,1,2)$, and a sheaf of principal ideals ${\mc I}_V$ in it such
that
$$
V=\Proj_{\POn}{\mc A}_V=\Proj_{\POn}{\mc O}_{\POn}[x,y,z,w]/{\mc I}_V
\subset\PA_{\POn}(1,1,1,2).
$$
Clearly, ${\mc I}_V$ is generated by a (weighted) homogeneous element
of degree 4, and $V\subset\PA_{\POn}(1,1,1,2)$ avoids a section along
which $\PA_{\POn}(1,1,1,2)$ is singular.

The second model is obtained by taking a double cover. We suppose
$Pic(V)=\ZA[-K_V]\oplus\ZA[F]$. Consider 
$$
{\mc E}=\rho_*{\mc O}_V(-K_V+mF),
$$
this is a vector bundle of rank 3 over $\POn$. We may choose $m$ such
that
$$
{\mc E}\simeq{\mc O}\oplus{\mc O}(n_1)\oplus{\mc O}(n_2),
$$
where $0\le n_1\le n_2$. Denote $b=n_1+n_2$.

We have the natural projection $\pi:X=\Proj_{\POn}{\mc E}\to\POn$. 
Let $M$ denote the class of the tautological bundle, 
$\pi_*{\mc O}(M)={\mc E}$, $L$ the class of a fiber of $\pi$, $l$ the
class of a line in a fiber, and $t_0$ the section corresponding to
the surjection ${\mc E}\to{\mc O}\to 0$, $t_0\circ M=0$. The following
diagram is commutative:
$$
\begin{CD}
V @>{\phi}>{R,2:1}> X\\
@V{\rho}VV @VV{\pi}V \\
\POn @= \POn
\end{CD}
$$
Here $\phi$ is a finite morphism of degree 2 branched over a divisor
$R\subset X$. $R$ is fibered into quartic curves over $\POn$, and we
may suppose
$$
R\sim 4M+2aL.
$$
The set $(n_1,n_2,a)$ defines $V/\POn$ up to moduli: different sets
corresponds to different varieties. While $n_1$ and $n_2$ can be
arbitrary chosen, $a$ should meet some lower estimation: $a$ can not
be much less than $n_1$, otherwise $R$ does not exist. 

Let us denote $H=\phi^*(M)$, $s_0=\frac12\phi^*(t_0$), 
$f=\frac12\phi^*(l)$; clearly, $F=\phi^*(L)$. As in the case $d=1$,
the divisor
$$
G_V\sim H-n_2F=(-K_V)+(a+n_1-2)F
$$
plays an important role in geometry of $V$. This is a minimally
twisted fibration into curves of genus 1 in $V$. Then, it is clear
that
$$
\overline{\bf NE}(V)=\RA_+[s_0]\oplus\RA_+[f],
$$
but now $s_0$ may not correspond to an effective 1-cycle on $V$. The
table below summarize the essential formulae on $V$ (see \cite{Grin1},
section 3.1):

\medskip
{\begin{center}
\renewcommand{\arraystretch}{1.5}
\begin{tabular}{|c|}
\hline
$K_V=-H+(a+b-2)F$ \\
\hline
$K_V^2=2s_0+(8-4a-2b)f$ \\
\hline
$(-K_V)^3=12-6a-4b$ \\
\hline
\end{tabular}
\end{center}}

\begin{remark}
We see again that all arguments and formulae work in the Gorenstein
case if we keep the assumption of non-reducity for fibers, so we may omit 
the non-singularity assumption. The unique difference with the case
$d=1$ is that Gorenstein singular $V/\POn$ may have reducible fibers.
Such fibers arise from twice covering of a plane with a double conic
as the ramification divisor, so the both described models are still 
suitable.
\end{remark}

\section{Fiber-to-fiber transformations}
\label{sec4}

This section is devoted to describing some special transformations of
del Pezzo fibrations. The birational classification problem formulated
in the bottom of section \ref{sec1}, is related to finding the set of
all Mori structures. In other words, we study birational maps modulo
birational transformations over the base. From the viewpoint of the
Sarkisov program, we only need links of type I, III, and IV, i.e.,
only that do change the Mori structure. So links of type II, which
does not change the structure of a fibration, are out of view at least
for the first look. Nevertheless, it is the type of the Sarkisov
links that is especially important in majority of questions related to
the birational classification problem, even if we only interested in
studying different Mori structures. The point is that proceeding with
the maximal singularities method or, more general, the Sarkisov
program and jumping from one Mori structure to another, we have to
know when we should stop. So, what we really need is a criterion or at
least strong enough sufficient conditions for a birational map to be
over the base, i.e., fiber-to-fiber. These conditions are known (see
\cite{Corti1}, theorem 4.2(i)), but the really hard thing is to prove
they are satisfied. Readers who handled any of various papers on the
maximal singularities method, may remember words like "excluding 
infinitely near singularities", this is exactly about the problem.
Birational rigidity and birational automorphisms are important
particular questions that are directly related to fiber-to-fiber
transformations. Thus, it is worth to know as much as possible about
this type of transformations, and this section gives some information
about that.

Let $V/C$ be a Mori fibration on del Pezzo surfaces of degree 1 or
2 over a curve $C$, and $\chi:V\dasharrow U$ a birational map onto 
another Mori fibration $U/C$ that is birational over the base, i.e.,
we have the following commutative diagram:
$$
\begin{array}{ccc}
V & \stackrel{\chi}{\dasharrow} & U \\
\downarrow && \downarrow \\
C & \stackrel{\simeq}{\lra} & C
\end{array}
$$

Taking the specialization at the generic point $\eta$ of $C$, we see
that $\chi$ induces a birational map
$$
\chi_{\eta}: V_{\eta}\dasharrow U_{\eta}.
$$ 
As it follows from \cite{Isk2}, theorem 2.6, $\chi_{\eta}$ is an
isomorphism if $d=K_{V_{\eta}}^2=1$, or decomposed into the so-called
Bertini involutions if $d=2$. Let us consider the last case in detail.

Bertini involutions can be constructed as follows. Let $A\in V_{\eta}$
be a rational point, i.e., defined over the field of functions of $C$.
There exists a morphism $\phi$ of degree 2 onto a plane (proposition
\ref{dp12surf_prop}) branched over a (non-singular) quartic curve.
Clearly $\phi$ defines an involution $\tau\in Aut(V)$ of $V_{\eta}$ that
transposes the sheets of the cover. Suppose $\phi(A)$ does not lie on
the ramification divisor, then there exists a point $A^*\in V_{\eta}$
which is conjugated to $A$ by means of $\tau$: $\tau(A)=A^*$. Clearly,
$A\ne A^*$. Take a pencil of lines on the plane through $\phi(A)$. 
On $V_{\eta}$ it becomes a pencil of elliptic curves (with degenerations) 
with exactly two base points: $A$ and $A^*$. Let us blow up these
points with exceptional divisors $e$ and $e^*$ lying over $A$ and $A^*$
respectively. We obtain an elliptic surface $S$ with 2 distinguished
sections $e$ and $e^*$. So, we have two (biregular) involutions
$\mu_A$ and $\mu_{A^*}$ on $S$ defined by these sections. Indeed, the
specialization of $S$ at the generic point of $e$ (or $e^*$, which is the
same) gives us an elliptic curve with 2 points $O$ and $O^*$ that 
corresponds to $e$ and $e^*$. Each of this points can be viewed as the
zero element from the viewpoint of the group law on elliptic curves,
thus we get two reflection $\mu'_A$ and $\mu'_{A^*}$ defined
as follows: for any point $B$
$$
\begin{array}{c}
B+\mu'_A(B)\sim 2O^*, \\
B+\mu'_{A^*}(B)\sim 2O.
\end{array}
$$
(just in this order: $\mu'_A$ is related to $O^*$, $\mu'_{A^*}$ to
$O$). These reflections give the biregular involutions $\mu_A$ and
$\mu_{A^*}$ on $S$ (note that $S$ is relatively minimal, so any fiber-wise
birational map is actually biregular). Finally, we blow down $e$ and $e^*$, 
and then $\mu_A$ and $\mu_{A^*}$ become the desired Bertini
involutions. Another way to describe $\mu_A$ is the following. Blow up
$A^*$, and you get a del Pezzo surface of degree 1. Its natural
involution defined by the double cover of a cone (see proposition
\ref{dp12surf_prop}) becomes birational on $V_{\eta}$, and this is
just our $\mu_A$.

Thus, referring to the case $d=2$, theorem 2.6 of \cite{Isk2} says:
$V_{\eta}$ is biregularly isomorphic to $U_{\eta}$, so we can view
$\chi_{\eta}$ as a birational automorphism of $V_{\eta}$, and then
there exists a finite set of points $I=\{A_1, A_2,\ldots,A_n\}$ on 
$V_{\eta}$ such that 
$$
\chi_{\eta}=\mu_{A_1}\circ\mu_{A_2}\circ\ldots\circ\mu_{A_n}\circ\psi,
$$
where $\psi\in Aut(V_{\eta})$. Now $\psi$ defines a fiber-wise 
birational transformation of $V/\POn$, hence applying it to the situation,
we get the case similar to degree 1, i.e., when we have an isomorphism
of the fibers over the generic point.

Thus, whether $d=1$ or $d=2$, the picture can be reduced to the
following case: we have the commutative diagram
\begin{equation}
\label{ftf_diag}
\begin{array}{ccc}
V & \stackrel{\chi}{\dasharrow} & U \\
\downarrow && \downarrow \\
C & \stackrel{\simeq}{\lra} & C
\end{array}
\end{equation}
where $\chi_{\eta}$ is an isomorphism of $V_{\eta}$ and $U_{\eta}$:
$$
\chi_{\eta}: V_{\eta}\stackrel{\simeq}{\lra} U_{\eta}.
$$
Now it is clear that if we throw out a finite number of points on $C$,
say, $P_1,P_2,\ldots,P_k$, then $\chi$ gives an isomorphism of $V$ and
$U$ over $C\setminus\{P_1,P_2,\ldots,P_k\}$. In its turn, $\chi$ can
be decomposed as
$$
\begin{array}{ccccccccccc}
V & \stackrel{\chi_1}{\dasharrow} & V_1 &
\stackrel{\chi_2}{\dasharrow} & V_2 & \stackrel{\chi_3}{\dasharrow} & 
\cdots & \stackrel{\chi_{k-1}}{\dasharrow} & V_{k-1} &
\stackrel{\chi_k}{\dasharrow} & U \\
\downarrow && \downarrow && \downarrow &&&& \downarrow && \downarrow \\
C & \stackrel{\simeq}{\lra} & C & \stackrel{\simeq}{\lra} & 
C & \stackrel{\simeq}{\lra} & \cdots & \stackrel{\simeq}{\lra} & 
C & \stackrel{\simeq}{\lra} & C 
\end{array}
$$
where $\chi_i$ is an isomorphisms over $C\setminus\{P_i\}$. In order
to distinguish such birational transformations from transformations
like Bertini involution, we will call them {\it fiber
transformations}. We may be naturally asked whether they are possible.
In order to clarify the question, let us consider 2-dimensional
situation. First, for any ruled surface, as an example of a fiber
transformation we can take any elementary transformation of a fiber:
blow up a point in a fiber and then contract the strict transform of
the fiber, which is a $(-1)$-curve. Note that the exceptional
divisor takes place of the original fiber. Obviously, such a
transformation is not an isomorphism, though it gives an isomorphism
of the fibers over the generic point of the base. On the other hand,
consider any two relatively minimal elliptic surfaces that are
birational over the base. It is very well known that such a birational
map is actually an isomorphism. This is true for any relatively
minimal (i.e., there are no $(-1)$-curves in fibers) fibrations into
curves of positive genus.

What about del Pezzo fibrations? On the one hand, they are fibered 
into rational surfaces, and we may expect they behave like ruled
surfaces. On the other hand, as it follows from results of the next
two sections, their essential birational properties are far away
from rational varieties, moreover, in many senses they behave like
elliptic fibrations. In fact, the very situation is a bit more 
complicated.

The rest of this section is contained in \cite{Grin3}, section 4 (as
to the case $d=1$), or \cite{Grin4} (the case $d=2$). Let us first
consider the case $d=1$.

So, we have the commutative diagram (\ref{ftf_diag}), $\chi$ is an
isomorphism of the fibers over the generic points. The situation can
be easy reduced to the following case. We assume $C$ to be a germ of a
curve with the central point $O$, $V$ and $U$ Gorenstein relatively
projective varieties over $C$ fibered into del Pezzo surfaces of
degree 1, their central fibers (i.e., over the point $O$) $V_0$ and
$U_0$ reduced. Algebraically, all that means the following. Let 
${\mc O}$ be a DVR (discrete valuation ring) with the maximal ideal
${\mathfrak m}=(t){\mc O}$, where $t$ is a local parameter, $V$ and
$U$ non-singular del Pezzo surfaces of degree 1 defined over 
${\mc O}$. If needed, we may take the completion of ${\mc O}$. Then,
let $K$ be the field of quotients, $C=\Spec {\mc O}$, $\eta=\Spec K$
the generic point of $C$. Using proposition \ref{dp12surf_prop}, we
may suppose that $V$ and $U$ are embedded into two copies $P$ and
$R$ of the weighted projective space $\PA_{\mc O}(1,1,2,3)$ respectively.
Denote $[x,y,z,w]$ and $[p,q,r,s]$ the coordinates in $P$ and $R$
with the weights (1,1,2,3).

By the condition, $\chi_{\eta}:V_{\eta}\to U_{\eta}$ is an
isomorphism. The key point is that $\chi_{\eta}$ induces an
isomorphism between $P_{\eta}$ and $R_{\eta}$. It follows easy from
the Kodaira vanishing theorem and the exact sequence of restriction
for ideals that define $V_{\eta}$ and $U_{\eta}$ as surfaces in
$P_{\eta}$ and $R_{\eta}$ respectively (see subsection 4.1,
\cite{Grin3}).

Then, we may choose the coordinates in $P$ and $R$ such that
\begin{equation}
\label{VU_eq}
\begin{array}{c}
V=\{w^2+z^3+zf_4(x,y)+f_6(x,y)=0\}\subset P, \\
U=\{s^2+r^3+rg_4(p,q)+g_6(p,q)=0\}\subset R,
\end{array}
\end{equation}
where $f_i$ and $g_i$ are homogeneous polynomials of degree $i$. Since
$\chi_{\eta}$ is an isomorphism of $P_{\eta}$ and $R_{\eta}$, it is
easy seen that $\chi$ and $\chi^{-1}$ can be defined as follows:
$$
\begin{array}{cc}
\chi=
\left\{
\begin{array}{ccl}
p & = & t^ax \\
q & = & t^by \\
r & = & t^cz \\
s & = & t^dw
\end{array}
\right\},
&
\chi^{-1}=
\left\{
\begin{array}{ccl}
x & = & t^{\alpha}p \\
y & = & t^{\beta}q \\
z & = & t^{\gamma}r \\
w & = & t^{\delta}s
\end{array}
\right\},
\end{array}
$$
where each of the sets $(a,b,c,d)$ and $(\alpha,\beta,\gamma,\delta)$
contains at least one zero. Then, all these numbers have to respect 
the graduation of $P$ and $R$, and, moreover, we know that $V$ and $U$ 
avoid singular points of $P$ and $R$, so for some integer $m>0$ the 
following conditions are satisfied:
$$
\left\{
\begin{array}{rcl}
a+\alpha & = & m  \\
b+\beta  & = & m  \\
c+\gamma & = & 2m \\
d+\delta & = & 3m \\
2d       & = & 3c \\
2\delta  & = & 3\gamma
\end{array}
\right.
$$
Using the symmetry of the situation, we may assume that $c=2k$,
$d=3k$, $\gamma=2l$, $\delta=3l$ with the conditions $k+l=m$ and 
$k\le l$ (i.e., $k\le\frac12m$). Now substituting these relations for
the numbers in (\ref{VU_eq}), we obtain
\begin{equation}
\label{fg_eq}
\begin{array}{c}
f_4(x,y)=t^{-4k}g_4(t^ax,t^by), \\
f_6(x,y)=t^{-6k}g_6(t^ax,t^by).
\end{array}
\end{equation}
Suppose $k=0$. Since the set $(a,b,c,d)$ is not composed from zeros
only (otherwise $\chi$ is already an isomorphism), then one of the 
numbers $a$ and $b$ must be positive. Let it be $a$. Then
(\ref{fg_eq}) shows that the equation
$$
w^2+z^3+zf_4(x,y)+f_6(x,y)=0
$$
defines a singularity of $V$ at the point $(t,x,y,z)=(0,1,0,0,0)$:
take an affine piece $x\ne 0$ (i.e., divide the equation by $x^6$) and
check that the differentials vanish at the indicated point.

Now let $k>0$. Since the set $(a,b,c,d)$ have to contain at least one
zero, we may assume $a=0$. Then $\alpha=m$. Since $l\ge k>0$,
$\alpha=m>0$, and the set $(\alpha,\beta,\gamma,\delta)$ contains
zero, we must suppose $\beta=0$. So $b=m-\beta=m$. Thus
$$
\begin{array}{c}
f_4(x,y)=t^{-4k}g_4(x,t^my), \\
f_6(x,y)=t^{-6k}g_6(x,t^my).
\end{array}
$$
It only remains to take into account that $k\le\frac12m$ and
$f_i\in{\mc O}[x,y,z,w]$, and we see that the equation for $V$ defines
a singular point at 
$$
(t,x,y,z,w)=(0,0,1,0,0).
$$ 
So we get the following assertion: $V$ must be singular if $\chi$ is 
not an isomorphism!

As to the case $d=2$ (\cite{Grin4}), we can repeat the previous arguments 
with the only difference that $P$ and $R$ have the type 
$\PA_{\mc O}(1,1,1,2)$, and $V$ and $U$ are defined by
$$
\begin{array}{c}
V=\{w^2+f_4(x,y,z)=0\}\subset P, \\
U=\{s^2+g_4(p,q,r)=0\}\subset R.
\end{array}
$$
Then $\chi$ and $\chi^{-1}$ have the form
$$
\begin{array}{cc}
\chi=
\left\{
\begin{array}{ccl}
p & = & t^ax \\
q & = & t^by \\
r & = & t^cz \\
s & = & t^dw
\end{array}
\right\},
&
\chi^{-1}=
\left\{
\begin{array}{ccl}
x & = & t^{\alpha}p \\
y & = & t^{\beta}q \\
z & = & t^{\gamma}r \\
w & = & t^{\delta}s
\end{array}
\right\},
\end{array}
$$
with the relation
$$
m=a+\alpha=b+\beta=c+\gamma=\frac12(d+\delta)
$$
for some $m>0$. Again, $(a,b,c,d)$ and $(\alpha,\beta,\gamma,\delta)$
are not consisted of only zeros, so by the symmetry we may assume 
$\gamma=0$ (thus $c=m$) and $d\ge\frac12m$. Finally, the relation
$$
f_4(x,y,z)=t^{-2d}g_4(t^ax,t^by,t^mz)\in{\mc O}[x,y,z]
$$
shows that the equation
$$
w^2+f_4(x,y,z)=0
$$
defines a singular point at
$$
(t,x,y,z,w)=(0,0,0,1,0),
$$
so $V$ is always singular if $\chi$ is not an isomorphism.

These results can be summarized as follows:
\begin{theorem}[Uniqueness of a smooth model]
\label{unique_th}
Let $V/C$ and $U/C$ be non-singular Mori fibrations on del Pezzo
surfaces of degree 1 or 2 that are birational over the base by means
of $\chi$:
$$
\begin{array}{ccc}
V & \stackrel{\chi}{\dasharrow} & U \\
\downarrow && \downarrow \\
C & \stackrel{\simeq}{\lra} & C
\end{array}
$$
Then $V$ and $U$ are always isomorphic, and $\chi$ is an isomorphism
if $d=1$, or a composition of Bertini's involutions if $d=2$. In other
words, any fiber transformation is trivial in this case. 
\end{theorem}

\begin{remark}
The similar result was proved for $d\le 4$ in \cite{Park} using
Shokurov's complements and connectedness principle. The difference is
that Park's result requires $V$ and $U$ to have non-degenerated
central fibers. In our case, we assume nothing about them.
\end{remark}

So, the reader can see that though $V/C$ is fibered into rational
surfaces, it behaves like elliptic fibrations in dimension 2: there
are no non-trivial fiber transformations without loss of smoothness.
Now it is time to give some examples.

\begin{example}["smooth case"] 
In these examples $U$ is non-singular.

First let $d=1$. Suppose $(a,b,c,d)=(0,6,2,3)$ and
$(\alpha,\beta,\gamma,\delta)=(6,0,10,15)$, $V$ and $U$ are defined by
$$
\begin{array}{cc}
V: w^2+z^3+x^5y+t^{24}xy^5=0, & U: s^2+r^3+p^5q+pq^5=0.
\end{array}
$$
It is easy to check that $U$ is non-singular, $V$ has a singular point
of type $cE_8$ (the so-called compound $E_8$-singularity) at
$(t,x,y,z,w)=(0,0,1,0,0)$. Note that the central fiber $V_0$ has a
unique singular point of type $E_8$, $U_0$ is non-singular.

The case $d=2$: suppose $(a,b,c,d)=(1,4,0,2)$,
$(\alpha,\beta,\gamma,\delta)=(3,0,4,6)$, $V$ and $U$ are given by
$$
\begin{array}{cc}
V: w^2+yz^3+tx^4+t^{12}y^4=0, & U: s^2+qr^3+tp^4+q^4=0.
\end{array}
$$
Again, $U$ is non-singular, $V$ has $cE_8$-singularity at
$(t,x,y,z,w)=(0,0,1,0,0)$. The central fiber $V_0$ is a non-normal del
Pezzo surface (the double cover of a cone branched over a triple plane
section), $U_0$ has an elliptic singularity.
\end{example}

\begin{example}["birational automorphism"]
The case $d=1$: suppose $(a,b,c,d)=(1,0,2,3)$,
$(\alpha,\beta,\gamma,\delta)=(0,2,2,3)$, and
$$
\begin{array}{cc}
V: w^2+z^3+t^4x^5y+xy^5=0, & U: s^2+r^3+p^5q+t^4pq^5=0.
\end{array}
$$
$V$ and $U$ have $cE_8$ singularity in the central fibers. Note that
$V$ and $U$ are biregularly isomorphic: put $w=s$, $z=r$, $x=q$, and
$y=p$. So we can assume $\chi$ to be defined as
follows:
$$
\begin{array}{cccc}
x\to t^{-1}y, & y\to tx, & z\to z, & w\to w.
\end{array}
$$
Thus, $\chi$ gives an example of a fiber transformation that is a
birational automorphism.

The degree 2: assume $(a,b,c,d)=(1,2,0,2)$,
$(\alpha,\beta,\gamma,\delta)=(1,0,2,2)$, $V$ and $U$ are given by
$$
\begin{array}{cc}
V: w^2+t^2y^3z+yz^3+x^4=0, & U: s^2+q^3r+t^2qr^3+p^4=0.
\end{array}
$$
$V$ and $U$ have $cD_4$-singularity. As before, they become isomorphic 
if you put $p=x$, $q=y$, $r=z$, $s=w$. So $\chi$ is a fiber
transformation that is a birational automorphism.
\end{example}

\section{Mori structures on del Pezzo fibrations: the case $d=1$}
\label{sec5}

In this section we formulate known results on the rigidity of
fibrations on del Pezzo surfaces of degree 1 and describe Mori
structures for non-rigid cases. We use the projective model via twice
covering described in subsection \ref{subsec32}.

\begin{theorem}[\cite{Grin1}, theorem 2.6]
Let $V/\POn$ be a non-singular Mori fibration on del Pezzo surfaces of
degree 1. Then $V/\POn$ is birationally rigid except for two cases:
\begin{itemize}
\item[1)] $\eps=0$, $n_1=n_2=n_3=2$;
\item[2)] $\eps=0$, $n_1=0$, $n_2=1$, $n_3=2$.
\end{itemize}
In other words, conjecture \ref{main_conj} holds for $d=1$.
Moreover, if $V/\POn$ is rigid, then this is a unique non-singular 
Mori fibration in its class of birational equivalency, as it follows
from theorem \ref{unique_th}.
\end{theorem}

\begin{remark}
It easy to check that linear systems $|n(-K_V)-F|$ are non-empty and
free from base components exactly for the listed cases. Indeed, note
that $|n(-K_V)-F|$ is non-empty and base components free if and only
if $|-K_V-F|$ has the same property (take the direct images of
anticanonical linear systems), and then use the projective model.
The case when $V/\POn$ met the $K^2$-condition was first proved 
in \cite{Pukh3}.
\end{remark}

So, the only two cases are non-rigid. Let us consider them in detail.

\subsection{The case $(\eps,n_1,n_2,n_3)=(0,2,2,2)$} First, let us
note that the distinguished section $s_b$ has the class $s_0$, and
this is a unique section with such a class of equivalency. Then, 
${\mc N}_{s_b|V}\sim{\mc O}(-1)\oplus{\mc O}(-1)$, so at least locally
there exists a flop centered at $s_b$. In fact, this flop leaves us
in the projective category: it is enough to check that the linear
system $n(-K_V)$ gives a birational morphism that contracts exactly
$s_b$. Thus, let $\psi:V\dasharrow U$ be such a flop. Consider the
linear system ${\mc D}=|-K_V-F|$ and its strict transform 
${\mc D}_U=\psi_*^{-1}{\mc D}$ on $U$. It is easy to see that 
$\Bas{\mc D}=s_b$ and a general member of ${\mc D}$ is nothing but a
del Pezzo surface of degree 1 that is blown up at the base point of
the anticanonical linear system. Moreover, $\dim{\mc D}=1$,
$\Bas{\mc D}_U=\emptyset$, and ${\mc D}_U$ is a pencil of del Pezzo
surfaces of degree 1. Thus, $U$ is fibered over $\POn$ onto del Pezzo
surfaces of degree 1, so
$$
\begin{array}{ccc}
V & \stackrel{\psi}{\dasharrow} & U \\
\downarrow && \downarrow \\
\POn && \POn
\end{array}
$$
is a link of type IV. It only remains to compute that the projective
model of $U/\POn$ (as a double cover) has the same structure constants
$(\eps,n_1,n_2,n_3)$ as $V/\POn$. We are ready to formulate the
following result:
\begin{proposition}[\cite{Grin1}, proposition 2.12]
Let $V/\POn$ be a non-singular Mori fibration on del Pezzo surfaces of
degree 1 with the set of structure constants
$(\eps,n_1,n_2,n_3)=(0,2,2,2)$. Then any other Mori fibration that is
birational to $V$, is birational over the base either to $V/\POn$
itself, or to $U/\POn$. $V/\POn$ and $U/\POn$ are the only
non-singular Mori fibrations in their class of birational equivalency.
In general case, $Bir(V)=Aut(V)\simeq\ZA_2$ and generated by an
automorphism corresponding to the double cover.
\end{proposition}

\begin{remark}
It may happen that $V$ and $U$ are isomorphic to each other (note they
have the same set of structure constants). In this case $\psi\in
Bir(V)$ and $Bir(V)\simeq\ZA_2\oplus\ZA_2$.
\end{remark}

\subsection{The case $(\eps,n_1,n_2,n_3)=(0,0,1,2)$, or a double cone
over the Veronese surface} Let $T\subset\PF$ be the Veronese surface
(i.e., $\PTw$ embedded into $\PF$ by means of the complete linear
system of conics), $Q\subset\PS$ a cone over $T$ with the vertex $P$,
and $R\subset\PS$ a cubic hypersurface such that $P\not\in R$ and
$R_Q=R\cap Q$ is non-singular. Then there exists a degree 2 finite 
morphism $\mu:U\to Q$ branched over $R_Q$. The variety $U$ (the
so-called double cone over the Veronese surface) is a
well-known Fano variety of index 2 with $\rho(U)=1$ and $(-K_V)^3=8$
(see the classification of Fano 3-folds, e.g., \cite{IskPr}). Then,
$U$ contains a two-dimensional family ${\mc S}$ of elliptic curves 
parametrized by $T$. These curves lie over the rulings of $Q$. ${\mc S}$
has one-dimensional sub-family of degenerations consisting of
rational curves with either a node or a cusp. 

Let $l\in{\mc S}$ be non-singular, and $\psi_l: V_l\to U$ the blow-up of
$l$. Then it is easy to see that $V_l$ is a non-singular Mori
fibration on del Pezzo surfaces of degree 1 with 
$(\eps,n_1,n_2,n_3)=(0,0,1,2)$. Conversely, take a del Pezzo fibration
$V/\POn$ with this set of structure constants. The linear system
$|-K_V-2F|$ consists of one element $G_V$ which is the direct product
of $\POn$ and an elliptic curve (\cite{Grin1}, lemma 2.9). Then the
linear system $|3(-K_V)-3F|$ defines a birational morphism 
$\psi:V\to U$ that contracts $G_V$ along the rulings.

But what if $l\in{\mc S}$ is singular? First, suppose $l$ is a
rational curve with an ordinary double point $B$. Let $\psi_1:U_1\to
U$ be the blow-up of $B$ and $E_1\sim\PTw$ the exceptional divisor.
Note that $l^1$ (the strict transform of $l$) intersects $E_1$ at two
points, and denote $t_1$ the line on $E_1$ that contains these points.
Now, blow up $l^1$: $\psi_2:U_2\to U_1$. It is easy to check that the
strict transform $t_1^2$ of $t_1$ on $U_2$ has the normal bundle isomorphic 
to ${\mc O}(-1)\oplus{\mc O}(-1)$, and we may produce a flop 
$\psi_3:U_2\dasharrow U_3$ centered at $t_1^2$, without loss of 
projectivity. The strict transform $E_1^3$ of $E_1$ on $U_3$ becomes
isomorphic to a non-singular quadric surface with the normal bundle
${\mc O}(-1)$, so $E_1^3$ can be contracted to an ordinary double
point: $\psi_4:U_3\to U_4=V_l$. We have the birational map
$\psi_l=(\psi_4\circ\psi_3\circ\psi_2\circ\psi_1)^{-1}:V_l\dasharrow
U$, and the reader can easy check that $V_l$ is a Gorenstein Mori
fibration (over $\POn$) on del Pezzo surfaces of degree 1 with a unique 
ordinary double point, and it has the structure set
$(\eps,n_1,n_2,n_3)=(0,0,1,2)$ (note that the constructions in section
\ref{sec3} work in the Gorenstein case as well, see remark
\ref{Gorenst_rem1}).

The case when $l$ has a cusp at a point $B\in U$ is similar to the
previous but a bit more complicated. As before, let $\psi_1:U_1\to U$
be the blow-up of $B$. We see that $l^1$ becomes non-singular and 
tangent to $E_1$ at some point $B_1\in E_1$. The tangent direction to 
$l_1$ at $B_1$ defines a line $t_1\subset E_1$. Now take the blow-up
$\psi_2:U_2\to U_1$ of the curve $l^1$. The strict transform $E_1^2$
of $E_1$ becomes isomorphic to a quadratic cone that is blown up at a 
point outside of the vertex of the cone. Moreover, the strict
transform $t_1^2$ is exactly a unique $(-1)$-curve on $E_1^2$ and has
the normal bundle ${\mc O}(-1)\oplus{\mc O}(-1)$. All that can be
checked as follows: blow up the point $B_1$ with an exceptional divisor
$E_2$, then blow up the strict transform of $l_1$. The strict
transform of $E_2$ is isomorphic to ${\mathbb F}_1$-surface and can be
contracted along its ruling. After that you get $U_2$. So, now we
produce a flop $\psi_3:U_2\dasharrow U_3$ centered at $l_1^2$. The
strict transform $E_1^3$ of $E_1^2$ is isomorphic to a quadratic cone
with the normal bundle ${\mc O}(-1)$, hence it can be contracted
$\psi_4:U_3\to U_4=V_l$ to a double point locally defined by the equation
$x^2+y^2+z^2+w^3=0$. We see that $V_l$ is a Gorenstein Mori fibration
over $\POn$ on del Pezzo surfaces of degree 1 with a unique singular 
point, and moreover, $V_l$ has the same structure set   
$(\eps,n_1,n_2,n_3)=(0,0,1,2)$ as before. Note that fibers of
$V_l/\POn$ are arised from elements of the pencil $|\frac12(-K_V)-l|$ 
on $U$.

About twenty years ago, S.Khashin (\cite{Kh}) tried to show the
non-rationality of $U$ by proving the uniqueness of the structure of a
Fano variety on it. Unfortunately, his arguments contain mistakes, and
up to now we have no a reliable proof of the non-rationality of $U$.
Nevertheless, the following conjecture seems to be true:
\begin{conjecture}
$U$ has exactly the following Mori structures: $U$ itself, and
$V_l/\POn$ for all $l\in{\mc S}$ (thus, "two-dimensional family" of
different del Pezzo fibrations). In particular, $U$ is not rational
(it has no conic bundles).
\end{conjecture}
Thus, it only remains to prove this conjecture, and the birational
identification problem for non-singular Mori fibrations on del Pezzo
surfaces of degree 1 will be completed.

\section{Mori structures on del Pezzo fibrations: the case $d=2$}
\label{sec6}

We first formulate the result and then describe the non-rigid cases.
The structure constants $(a,n_1,n_2)$ are taken from subsection
\ref{subsec33}. Recall $b=n_1+n_2$.

\begin{theorem}[\cite{Grin1}, \cite{Grin2}]
Let $V/\POn$ be a non-singular Mori fibration on del Pezzo surfaces of
degree 2. Then $b+2a>0$, and if $b+2a>2$, $V/\POn$ is birationally
rigid.

Suppose $b+2a=2$. Then the only following cases are possible:
\begin{itemize}
\item[1)] $a=0$, $n_1=0$, $n_2=2$;
\item[2)] $a=-2$, $n_1=2$, $n_2=4$;
\item[3)] $a=-3$, $n_1=2$, $n_2=6$;
\item[4)] $a=1$, $n_1=0$, $n_2=0$;
\item[5)] $a=0$, $n_1=1$, $n_2=1$;
\item[6)] $a=-1$, $n_1=2$, $n_2=2$.
\end{itemize}
General varieties in the cases 1)--3) are rigid. The cases 4)--6) are
all non-rigid. 

Suppose $b+2a=1$, then the only following cases are possible:
\begin{itemize}
\item[7)] $a=0$, $n_1=0$, $n_2=1$;
\item[8)] $a=-1$, $n_1=1$, $n_2=2$.
\end{itemize}
Both of them are non-rigid.

So, under the assumption of generality for the cases 1)--3),
conjecture \ref{main_conj} holds for $d=2$, and if $V/\POn$ is rigid, 
then this is a unique non-singular Mori fibration in its class of 
birational equivalency (theorem \ref{unique_th}).
\end{theorem}

\begin{remark}
The $K^2$-condition corresponds to the case $b+2a\ge 4$, and the
result was first obtained by A.Pukhlikov (\cite{Pukh3}). The reader
has not to be confused with the generality assumptions in the cases
1)--3) (see \cite{Grin2}), they arised from some technical troubles of 
the maximal singularities method. The author believes they can be 
omitted finally. 
\end{remark}

\subsection{The case $(a,n_1,n_2)=(1,0,0)$} In this case
$X\simeq\POn\times\PTw$, so $V$ is nothing but a double cover of
$\POn\times\PTw$ branched over a divisor $R$ of bi-degree $(2,4)$.
Thus, $V$ is a del Pezzo fibration with respect to the projection onto
$\POn$. On the other hand, the projection onto $\PTw$ represents $V$
as a conic bundle. Indeed, a fiber of $V\to\PTw$ is a double cover
of a line branched over 2 points, hence either a conic, or a couple of
lines, or a double line. It is easy to see that the discriminant
curve of $V\to\PTw$ has the degree 8. So, $V$ has at least two
different Mori structures. Note that $\POn\leftarrow V\to\PTw$ is a
trivial example of type IV links.

\subsection{The case $(a,n_1,n_2)=(0,1,1)$, or a double quadratic cone}
In this case the linear system $|-2K_V|$ gives a small contraction onto 
the canonical model of $V$, which can be realized as double covering of 
a non-degenerated quadratic cone $Q\subset\PQ$ branched along a quartic 
section $R_Q$. It easy to see that $V$ has at most two curves of the
class $s_0$. If a curve of the class $s_0$ is unique on $V$
(this means that $R_Q$ contains the vertex of $Q$), then $s_0$ is
the so-called {\bf -2}-curve of the width 2 (in the notions of
\cite{Reid2}). 
Otherwise (if $R_Q$ does not pass through the vertex of $Q$), there are 
two curves of the class $s_0$, which are disjoint and have the normal
bundles ${\mc O}(-1)\oplus{\mc O}(-1)$. In both the cases we 
obtain another structure of a non-singular Mori fibration on del Pezzo 
surfaces of  degree 2 after making a flop centered at these curves.
Moreover, the second structure has the same set $(a,n_1,n_2)$.
Note that both the Mori structures arise from two families (actually,
two pencils) of planes on $Q$.

If $R_Q$ does not pass through the vertex of the cone, it was proved
that a general variety of this type has exactly two Mori structures
just described (\cite{Grin5}). Actually, it should be true always, not
only for general cases, and even if $R_Q$ did pass through the vertex
of the cone, but the proof of this fact is still waiting for its time. 

\subsection{The case $(a,n_1,n_2)=(-1,2,3)$}
First, let us note that $V$ has a unique curve of the class $s_0$
(say, $C$), and it can be contracted (small contraction) by the linear 
system $|nH|$ for $n\ge 2$. The normal bundle of this curve is isomorphic 
to  ${\mc O}(-1)\oplus{\mc O}(-2)$, hence there exists an anti-flip
$\psi:V\dasharrow U$ centered at $C$ (this is the simplest example of
anti-flip, the so-called Francia's anti-flip: blow up $C$,
make a flop centered at the minimal section of the
corresponding ruled surface, and then contract the strict transform of
the exceptional divisor, which is $\PTw$ with the normal bundle 
${\mc O}(-2)$). $U$ has a unique (non-Gorenstein) singular point of
index 2 (the latter means that $2K_U$ is a Cartier divisor). Moreover,
$U$ turns out to be a Mori fibration over $\POn$ on del Pezzo surfaces 
of degree 1. It is not very difficult to check that general elements of the
pencil $|-K_V-F|$ are del Pezzo surfaces of degree 1 that are blown up
at the base point of the anticanonical linear system, and $C$ is
exactly their common $(-1)$-curve. Thus, $\psi$ make them to be the 
fibers of $U/\POn$. Conjecturally, $V/\POn$ and $U/\POn$ are the only Mori
structures in this case.

\subsection{The case $(a,n_1,n_2)=(0,0,1)$, or a double space of index
2} This case is similar to a double cone over the Veronese surface.
Let $U$ be a double cover of $\PT$ branched over a smooth quartic.
This is a well known Fano variety of index 2, the so-called double
space of index 2. Denote $H_U$ the generator of the Pickard group,
$\Pic(U)=\ZA[H_U]$. Clearly, $K_U\sim -2H_U$. Let $l\subset U$ be 
a curve of genus 1 and degree 2 (i.e., $H\circ l=2$), so $l$ is a
double cover of a line in $\PT$ branched at four points. Then $l$ is
either a (non-singular) elliptic curve, or a rational curve with a
double point (node or cusp), or a couple of lines with two points of 
intersection (probably, coincident).

Suppose $l$ is non-singular, and $\psi_l:V_l\to U$ is the blow-up of
$l$. Then $V_l$ is a non-singular Mori fibration on del Pezzo surfaces
of degree 2. The fibers of $V_l/\POn$ arise from a pencil of planes in
$\PT$ that contain the image of $l$. It is not very difficult to
compute that $V_l/\POn$ has the structure set $(a,n_1,n_2)=(0,0,1)$.
Conversely, given $V/\POn$ with such a structure set, we can contract
a unique divisor in $|-K_V-2F|$, which is isomorphic to the direct
product of an elliptic curve and a line, and get a double space of
index 2.

If $l$ is singular but irreducible, i.e., a rational curve with a node
or a cusp, we can also obtain a Gorenstein Mori fibration on del Pezzo 
surfaces of degree 2, but now with a singular point of type
respectively $x^2+y^2+z^2+w^2=0$ or $x^2+y^2+z^2+w^3=0$, "blowing up"
the curve $l$ just in the way as we obtain fibered structures from a
double cone over the Veronese surface (see the previous section).
$V_l/\POn$ has the same set of structure constants:
$(a,n_1,n_2)=(0,0,1)$.

However, $U$ has another type of Mori fibration. Suppose $l$ is a
couple of lines (say, $l=l_1\cup l_2$) intersected by two different 
points. Blow up one of this line and make a flop centered at the strict 
transform of the second line, and you obtain a structure of a 
(non-singular) Mori fibration on cubic surfaces $V_{l_1}/\POn$. 
As before, its fibers arise from a pencil of planes that contain the 
image of any of these lines in $\PT$. The same construction works if 
$l$ consists of two lines that are tangent to each other. Note that we
can change the order of the lines (first blowing up $l_2$) and get
$V_{l_2}/\POn$. The point is that $V_{l_1}/\POn$ is (biregularly) 
isomorphic to $V_{l_2}/\POn$ over the base, and this order change
corresponds to a birational automorphism (the so-called Geiser
involution, see \cite{Isk2}). So we may denote these fibration on
cubic surfaces $V_l/\POn$ without any confusion.

Following to \cite{Kh2}, we can formulate the following conjecture:
\begin{conjecture}
Let $U$ be a double space of index 2. Then $U$ has the following Mori
structures: $U$ itself, and $V_l/\POn$ for all curves $l$ of degree 2
and genus 1.
\end{conjecture}

In this subsection it only remains to note that $U$ is known to be
non-rational, as it follows from \cite{Tikh}.

\subsection{The case $(a,n_1,n_2)=(-1,1,2)$, or a singular double cone
over the Veronese surface}
As to this case, note first that the linear system $|H-2F|$ has a
unique representative, which we denote $G_V$. Then, there is a unique 
curve of the class $s_0$ on $V$, and this curve has the normal bundle
${\mc O}(-1)\oplus{\mc O}(-1)$. So we produce a flop (actually, in the
projective category) $\psi_1:V\dasharrow V^+$. The strict transform
$G_V^+$ of $G_V$ on $V^+$ is isomorphic to either $\POn\times\POn$
or a quadratic cone, and has the normal bundle ${\mc O}(-1)$ in both
the cases. So there exist a contraction $\psi_2:V^+\to U$, which gives
us a Fano variety with a double point of type $x^2+y^2+z^2+w^2=0$ or 
$x^2+y^2+z^2+w^3=0$ respectively. It is easy to check that this $U$
can be obtained as a double cover of the cone over the Veronese
surface branched over a cubic section that does not pass through the
vertex of the cone and has a unique du Val point of type $A_1$ or
$A_2$. Conversely, given $U$ with such a kind of singularity, we
can always obtain a structure of del Pezzo fibration as indicated.
Now it is clear that $U$ can be transformed to a (singular Gorenstein)
Mori fibration on del Pezzo surfaces of degree 1 as in the previous
section. Conjecturally, these are all Mori structures on $U$ (i.e.,
one structure of Fano variety, one structure of fibration on del Pezzo
surfaces of degree 2, and "2-dimensional family" of structures of
fibrations on del Pezzo surfaces of degree 1), and hence $V$ is a
unique non-singular Mori fibration in its class of birational
equivalency.

\end{document}